\documentclass[12pt]{article}
 \usepackage{amsmath,amsfonts,theorem}
 \newcommand{\ntag}[1]{} 
\numberwithin{equation}{section}
 
 \newtheorem{prop}{Proposition}

 \newtheorem{cor}{Corollary}

 \theorembodyfont{\rmfamily}
 \newtheorem{dfn}{Definition}

 \newcommand{\qed}{\ifhmode\unskip\nobreak\fi\quad\ensuremath\square}

 \newcommand{\sA}{\mathcal A} 
 \newcommand{\sG}{\mathcal G} 
 
 \newcommand{\sI}{\mathcal I}

 \newcommand{\Oh}{\mathcal O}

 \newcommand{\ga}{\gamma}
 
 \newcommand{\om}{\omega}

 \newcommand{\La}{\Lambda}
 \newcommand{\Om}{\Omega}
 \newcommand{\la}{\lambda}

 \newcommand{\C}{\mathbb C}

 \newcommand{\R}{\mathbb R}
 \newcommand{\Z}{\mathbb Z}

\begin{document}

 \title{Universal Maslov class of Bohr - Sommerfeld lagrangian
 embedding to pseudo - Einstein manifold}

  \author{N.A. Tyurin}

\date{}

\maketitle

\begin{center}
 {\em BLThP JINR (Dubna) and MSU of RT,}\\
 {\it ntyurin@theor.jinr.ru}
 \end{center}
 \bigskip

\section*{Introduction}

In this paper we define some universal  1 - cohomology class on
a lagrangian submanifold $S \subset M$ of a simply connected
compact symplectic manifold $(M, \om)$ which satisfy the following
conditions:

--- $(M, \om)$ is pseudo - Einstein s.t.
$$
K_{\om} = k . [\om] \in H^2(M, \Z),
$$
where $K_{\om}$ is the associated canonical class of the symplectic form $\om$
and $k$ is a constant;

--- $S$ is Bohr - Sommerfeld (the definition of
Bohr - Sommerfeld lagrangian submanifolds see in Section 1 or in
\cite{GT}).

We call this 1 - cohomology class the universal Maslov class since, as
it is explained below, this class is a natural generalization
of the known Maslov class for lagrangian immersions to
symplectic vector spaces, see \cite{ArGiv}. At the same time our construction
doesn't follow  the standard  for general for the compact case way which uses  the notion of
 index for membranes which bound loops on a given lagrangian submanifold, introduced 
  in \cite{KarMas} and exploited in many other papers (we refer just to  one of them, [2],
  since the situation we study here is non integrable vesrion of the situation, considered
  there). The index of membranes was called the Maslov index and is exploited in
  pure symplectic case. On the other hand, 
 Fukaya in \cite{Fuk} already mentioned that in
the situation which we study it is possible to define some 
index for Bohr - Sommerfeld lagrangian submanifolds of a pseudo - Einstein
manifold with respect to an integrable complex structure and this index generalizes 
the notion of the Maslov class.
In our construction one doesn't require that such a structure exists on
our symplectic manifold $M$. This paper is 
a consequence of \cite{jT3} where one mostly  concentrates on the case
of Kahler - Einstein manifolds since in this case there is some relationship
of the constructed Maslov class and the minimality problem
for lagrangian submanifolds. At the same time it is clear (we will see it at Section 3)
that the definition of the universal Maslov class, given in the present paper,
agrees with the definition of the Maslov class of a Bohr - Sommerfeld lagrangian embedding
to a Kahler - Einstein manifold, given at \cite{jT3}. Namely, if a given simply
connected pseudo - Einstein symplectic manifold admits a Kahler - Einstein metric,
then for a Bohr - Sommerfeld lagrangian embedding the universal Maslov class
and the Maslov class, defined by the Kahler - Einstein metric, coincide.
It follows that we can add  several new facts to the discussion
of \cite{jT3}. First, for a given simply connected Kahler - Einstein
manifold the Maslov class of a Bohr - Sommerfeld lagrangian submanifold doesn't
depend on the choice of the Kahler - Einstein structure. Second, if a Bohr
- Sommerfeld lagrangian submanifold has non trivial universal Maslov class in
a simply connected  pseudo - Einstein symplectic manifold then
there is no Kahler - Einstein metric on this symplectic manifold such that
the lagrangian submanifold is minimal with respect to it.

In Section 4 we compare the universal Maslov class and the Maslov index.
They are different although they have similiar geometric interpretations.
And the difference is given by a multiple of the symplectic area which itself gives
an integer class on the group $\pi_2(M, S; \Z)$ if $S$ is Bohr - Sommerfeld. 
 
I would like to thank D. Orlov and D. Auroux for  helpful discussions and comments.
This work was done due to the partial financial support of RFBR (grants NN 05 - 01 -
01086, 05 - 01 - 00455).

\section{Bohr - Sommerfeld condition}

Let $(M, \om)$ be a simply connected compact symplectic manifold and suppose
that its canonical class $K_{\om}$ is proportional to the cohomology class
of the symplectic form:
$$
K_{\om} = k. [\om] \in H^2(M, \Z).
$$
We call such a manifold pseudo - Einstein, following \cite{Fuk},
although in some cases (if $k<0$)
 it is reasonable to use another term --- pseudo - Fano symplectic manifolds.

We will work with the anticanonical line bundle $K^{-1}_{\om} \to M$
denoting it as $K^{-1}$. It is realized as follows:
let us fix an almost complex structure $I$, compatible with $\om$,
then the pair $(\om, I)$ can be completed to the corresponding hermitian triple
$(\om, G, I)$, where $G$ is the corresponding riemannian metric.
Thus it induces a hermitian structure $H$ on the tangent bundle $TM$
and as a complex bundle $(TM, H)$ is isomorphic to the holomorphic
tangent bundle $T^{1,0}M$ (which is well defined for non - integrable complex
structures as well as for the integrable ones). The determinant complex line bundle
$\det T^{1,0} M$, endowed with the corresponding hermitian structure,
is the anticanonical line bundle $K^{-1}$. Topologically as a $\C^*$
- bundle it is defined by the first Chern class $c_1(K^{-1})$
which we denote by the same symbol $K^{-1}$, following  algebro - geometrical
traditions.

For any almost complex structure one has the corresponding hermitian structure
$h$ on the complex line bundle $K^{-1}$ and thus the space of hermitian connections
$\sA_h(K^{-1})$ is defined (all details of the theory of connections,
curvatures  and gauge transformations can be found, f.e., in \cite{DK}).
The gauge group $\sG_h$ acts on this space and it is known that for the abelian
connections on a simply connected manifold every gauge equivalence class of
connections is defined by the curvature
form which belongs to $\Om^2_M(i \R)$, see \cite{DK}. Therefore if we
impose the condition on the curvature form
$$
F_a = - 2 \pi i k \cdot \om,
\eqno (1)
$$
we get precisely one orbit
$$
\sG_h(a) = \Oh_{k \cdot \om} \subset \sA_h,
$$
where $a$ is a connection, satisfies (1), and this orbit consists of all solutions
of (1). As usual, for each pair $a, a_1 \in \Oh_{k \cdot \om}$ the
difference $a - a_1 \in \Om^1_M(i \R)$ is a pure imaginary
exact 1 - form.

Let us fix an element $a$ from $\Oh_{k \cdot \om}$. Then for
a smooth lagrangian submanifold $S \subset M$ one can consider
the restriction of the pair $(K^{-1}, a)$ on it getting
a trivial line bundle with a flat connection. Indeed,
since the curvature $F_a$ is proportional to the symplectic
form, $a$ is flat being restricted to any lagrangian submanifold
by the definition. Thus we get a pair $(K^{-1}|_S, a|_S)$ and
for this flat connection one considers its character
on the fundamental group $\pi_1(S)$.

We call here a lagrangian submanifold {\it Bohr - Sommerfeld with respect to
the anticanonical bundle}  (or just Bohr - Sommerfeld for short) if this character
is trivial (see, f.e., \cite{GT}, and the references therein). Usually in geometric quantization
one uses this notion with respect to the prequantization bundle (see, f.e., \cite{GT})
but being pure geometrical it can be used in more general setup with respect to any line
bundle whose first Chern class is proportional to the cohomology class of the symplectic
form. Let's remark that we do not require that the last one is integer since we are working
with the anticanonical bundle.  At the same time in the usual for geometric quantization
situation when the class of symplectic form is integer every standard Bohr - Sommerfeld
lagrangian submanifold is Bohr - Sommerfeld with respect to the anticanonical class,
therefore the using of the same term is natural.

At the first glance the definition depends on the choice
of the reference connection $a \in \sA_h(K^{-1})$ but
the point is that for a simply connected pseudo - Einstein symplectic
manifold this notion is universal.

Indeed, if we change the reference connection $a$ in the gauge equivalence class
$\Oh_{k \cdot \om}$ then the character must be the same for any $S \subset M$ since it is
invariant under the gauge group action. If, further, we change the hermitian structure
$h$ to another $h'$ on the anticanonical line bundle, then connections $a \in \Oh_{k \cdot \om}
\subset \sA_h(K^{-1})$ and $a' \in \Oh_{k \cdot \om}' \subset
\sA_{h'} (K^{-1})$ belong to the same $\C^*$ - equivalence class
in the big space $\sA(K^{-1})$ of all $\C^*$ - connections on $K^{-1}$
 since they have the same curvature form and thus the characters of $a|_S$
 and $a'|_S$ must be the same.

 Thus in the situation which is studied   the present paper the notion of
 Bohr - Sommerfeld with respect to the anitcanonical  bundle
 lagrangian submanifolds is universal.

 It was discussed several times, see \cite{jT1}, \cite{jT2},
 that while the standard lagrangian condition is static,
 the standard  Bohr - Sommerfeld condition is dynamical and the same is true for
 the Bohr - Sommerfeld condition with respect to the anticanonical bundle; this means
 that locally the space of Bohr - Sommerfeld lagrangian submanifolds
 is generated by strictly hamiltonian vector fields
 and therefore the space of all  isodrastic deformations of a given
 Bohr - Sommerfeld lagrangian submanifold is exhausted by infinitesimal
 symplectomorphisms of $(M, \om)$ (and of course every given $S_0$ in this
 representation has a huge stabilizer ${\rm Stab} S_0 \subset
 {\rm Sym} (M, \om)$).

In what follows we work with Bohr - Sommerfeld with respect to
the anticanonical bundle lagrangian submanifolds, omitting for breviety
the characterization.

\section{Universal class}

Let $(M, \om)$ be a simply connected
pseudo - Einstein manifold and $S$ is a
smooth lagrangian submanifold, satisfies the Bohr - Sommerfeld
condition with respect to the anticanonical line bundle
(note that every Bohr - Sommerfeld lagrangian submanifold must be
orientable).

Let us fix an almost complex structure $I$, compatible with $\om$,
and take an abelian connection $a \in \Oh_{k \cdot \om}
\subset \sA_h(K^{-1})$ as it was done in section 1. On the other hand,
the choice of $I$ gives an isomorphism
$$
K^{-1}|_S = \det TS \otimes \C
\eqno (2)
$$
together with an identification of the hermitian structure
on $K^{-1}|_S$ induced by $I$ and the complexification
of the special orthogonal structure on $\det TS$, given
by the restriction to $S$ of the riemannian metric $G$
from the hermitian triple $(\om, G, I)$.  It can be easily shown
(see, f.e., \cite{jT3}) that the canonical trivialization of the restriction
$K^{-1}|_S$ to an orientable lagrangian submanifold $S$ is given by the orthogonal
projection of the top polivector field of unit lenth to $\La^n T^{1,0} M|_S$. Indeed,
at any point $p \in S \subset M$ the space
$$
\La^n T^{\C}_p M
$$
contains

--- real line $\La^n T_p S$;

--- line without real points $\La^n T^{1, 0}_p M = K^{-1}|_p$;

--- the hermitian product $H$.

Local computations show that the orthogonal projection of a unit vector from
$\La^n T_p S$ to the complex line $K^{-1}_p$ never  vanishes. If $S$ is orientable
then $\La^n T_p S$ is trivial and there are two sections of unit lenth, whose orthogonal
projections to $K^{-1}|S$ induce two trivializations, which we call canonical,
and since these trivializations are conjugated by the canonical $U(1)$ - action,
we have one canonical trivial connection on $K^{-1}|_S$ which we denote as
$A_0 \in \sA_h(K^{-1}|_S)$ and which is independent on the orientation choice.
 Thus under the indentification (2) there are two flat connections
 with trivial characters: the restriction of the reference connection $a$ and the canonical connection $A_0$.  Therefore as in paper
\cite{jT3} we can compare two flat connections with trivial characters
on $\pi_1(S)$;
it follows from the coincidence of the characters that
$a|_S$ and $A_0$ belong to the same class
of  gauge equivalence in $\sA_h(K^{-1}|_S)$ and hence they
are related by a gauge transformation
$$
g(a|_S, A_0) = g_S \in {\rm Map} (S, U(1)).
$$
This gauge transformation gives us some 1 - cohomology class on $S$
by the rule
$$
H^1(S, \Z) \ni m_S = g_S^* h,
\eqno (3)
$$
where $h \in H^1(U(1), \Z)$ is the generator of $H^1(U(1), \Z)$.

From this we have
\begin{dfn} The universal Maslov class $m_S \in H^1(S, \Z)$, induced by a Bohr - Sommerfeld
lagrangian embedding $S \subset M$,
is given by the equation (3).
\end{dfn}

{\bf Example.} If one takes a symplectic vector space $V^{2n}$ with a constant symplectic
form $\Om$ as the simply connected symplectic manifold then the universal
Maslov class, introduced above, is exactly the "classical" Maslov class from \cite{ArGiv}.
Indeed, in this case one can take as the reference connection on the anticanonical line
bundle the determinant Levi - Civita connection of a constant integrable complex structure $I$,
compatible with $\Om$. Then the universal class (3) can be computed from the canonical pairing
of the global  covariantly constant section, trivializing  the anticanonical line bundle,
and the riemannian volume form on the lagrangian submanifold. And this  gives the classical
definition.

We must emphasize however that at the moment our universal Maslov class exists
only for smooth embeddings $S \subset M$ {\it by the definition}. Indeed, for the immersion
case one should work with singular connections and
gauge transformations (singularities arise at the self intersection points of the images
of immersions).

The class introduced above is universal since

\begin{prop} The definition of the universal Maslov class of a Bohr - Sommerfeld
lagrangian submanifold $S \subset M$ of a pseudo - Einstein
simply connected symplectic manifold $(M, \om)$ is correct,
s.t. it doesn't depend on the choice of $I$ and $a$.
\end{prop}

Indeed, the space of compatible almost complex structure for
a symplectic manifold is contractible. For each almost complex structure
giving the corresponding hermitian structure $h$ on $K^{-1}$ the orbit
$\Oh_{k \cdot \om}$ is connected (and the simply connectedness of
$M$ is essential!). The topological type of $g_S$ is constant on
the connected space of pairs $(I, a)$, and this type is
exactly carried by  the class $m_S$.

Moreover, it's not hard to establish that
\begin{prop} 1. The universal Maslov class of a Bohr -
Sommerfeld lagrangian submanifold of a
simply connected  pseudo - Einstein
 symplectic manifold is invariant under
isodrastic deformations.

2. The universal Maslov class of a Bohr - Sommerfeld
lagrangian submanifold of a simply connected pseudo - Einstein
symplectic manifold is equivariant with respect to the action of
the  symplectomorphism group
of $(M, \om)$.
\end{prop}

{\bf Remark 1.} It's not hard to extend the considerations to
non simply connected case. Indeed, if $M$ has non trivial 1- cohomology group:
$$
b^1(M) > 0,
$$
then there are various $\sG_h$ - orbits of solutions for equation (1) enumerated by
the real cohomology space $H^1(M, \R)$
of our given $M$ (if we choose some origin
in the affine space $\sA_h(K^{-1})$, associated with the vector space
$\Om^1_M$). For a point $[a] \in H^1(M, \R)$ we denote the corresponding
orbit as $\Oh_{k \cdot \om}^{[a]}$. Then if an orientable lagrangian submanifold
$S \subset M$ is Bohr - Sommerfeld with respect to a class $[a]$ then we can
construct some class $m_S = m_S(a)$ using the same strategy as it was in the
simply connected case.
But in this case it can happen that $S \subset M$ is Bohr - Sommerfeld for
different connection classes; for a pair of such classes $[a_1], [a_2]$, represented
by abelian connections $a_1, a_2 \in \sA_h(K^{-1})$ we have that their difference
$a_1 - a_2 \in \Om^1_M(i \R)$ is a closed 1 - form. After restriction
to $S$ this form has integer values on $H_1(S, \Z)$ since $S$ is Bohr - Sommerfeld
for both $a_1$ and $a_2$. Therefore it defines an element from $H^1(S, \Z)$ and it's
clear that this element equals exactly to the difference
$$
m_S(a_2) - m_S (a_1) \in H^1(S, \Z).
$$
Hence the class, given by our construction, is not universal
in general for non simply connected case; however if the evaluation map
$$
H^1(M, \R) \times H_1(S, \Z) \to \R
$$
is trivial, the construction works and we get the universal class. Otherwise we have
some interesting "tunneling effect".

{\bf Remark 2.} The universal Maslov class can be understood as an obstruction.
Consider a simply connected pseudo - Einstein symplectic manifold $(M, \om)$ together with
a Bohr - Sommerfeld lagrangian submanifold. If we fix a compatible almost complex structure
then we get the canonical trivial connection $A_0 \in \sA_h(K^{-1}|_S)$. The natural question
arises: can this connection be extended to a global connection $a \in \sA_h(K^{-1})$
with the curvature form, proportional to the symplectic form? If the universal
Maslov class $m_S \in H^1(S, \Z)$ is non trivial then the extension
doesn't exist. On the other hand, if the universal Maslov class is trivial
$$
m_S = 0,
$$
then such an extension exists (and there are lot of such extensions). Thus the universal
Maslov class is the obstruction to the existence of such extensions.

There are some additional remarks in the non orientable case, but we leave it outside
of our present discussion.

\section{The Kahler - Einstein case}

Consider now the following situation: let $(M, \om)$ be a simply connected
pseudo - Einstein manifold, $S \subset M$ be a Bohr - Sommerfeld lagrangian
submanifold and suppose that $(M, \om)$ admits Kahler - Einstein metrics.
This means that there are exist integrable complex structures compatible with
our symplectic form $\om$ and moreover the corresponding Kahler metrics
have the same Ricci form, proportional to $\om$, see \cite{GrHar}.
The case of lagrangian embeddings to Kahler manifolds was studied in
\cite{jT3}, and there one establishes that the classical definition
of the Maslov class from \cite{ArGiv} can be generalized to the case
of lagrangian embeddings satisfy some strong property:
the restriction to $S \subset M$ of the determinant Levi - Civita connection
$a_{LC}$ is flat and trivial (admits covariantly constant sections), see
\cite{jT3}. For such a lagrangian submanifold one defines the phase
$$
g_S: S \to U(1)
$$
which is the gauge transformation, transporting $a_{LC}|_S$ to $A_0$, and
the Maslov class is given by the formula
$$
H^1(S, \Z) \ni m_S = g_S^* h, h \in H^1(U(1), \Z),
$$
details see in \cite{jT3}.

Denote as $\sI$ the moduli space
of Kahler - Einstein metrics on $M$ and consider
an element $I$ of this space. In our case for a Bohr - Sommerfeld lagrangian
submanifold $S \subset M$ the flatness and the triviality conditions
on the determinant Levi - Civita connection $a_{LC}$, induced by $I$,
are satisfied automatically and it is not hard to see that

\begin{prop} The Maslov class $m_S = m_S(I)$ induced by the complex structure $I$
coincides with the universal Maslov class.
\end{prop}

Indeed, for a Kahler - Einstein metric the determinant Levi - Civita connection
$a_{LC} \in \sA_h(K^{-1})$ belongs to the orbit $\Oh_{k \cdot \om}$ defined in Section 1,
and thus $a_{LC}$ can be taken as the reference connection $a$ of the construction
of Section 2. Then the 1- cohomology class is the same by the definition.

But the universal Maslov class doesn't depend on the choice of complex structures and it follows
that the Maslov class of \cite{jT3} is the same for all elements of $\sI$. Further
we've established in \cite{jT3} that  the Maslov class in the Kahler  - Einstein
case is the obstruction
to the possibility of isodrastic deformation  to a minimal lagrangian submanifold.
As a corollary we get the following

\begin{prop} Let $(M, \om)$ be a simply connected pseudo - Einstein symplectic manifold
and $S \subset M$  be a Bohr - Sommerfeld lagrangian submanifold. Then if the universal
Maslov class $m_S \in H^1(S, \Z)$ is non trivial, then there is no Kahler - Einstein metric
on $M$ such that $S$ is minimal with respect to it.
\end{prop}

The proof is obvious.

\section{Universal Maslov class, Maslov index and monotone lagrangian submanifolds}

In general situation for a lagrangian submanifold in a symplectic manifold one can define some
index for loops which depends on the elements of the group $\pi_2(M, S; \Z)$. 
And it is natural to compare the Maslov index for loops with the universal Maslov class
defined in the case of Bohr - Sommerfeld embeddings to pseudo - Einstein symplectic manifolds.

Let $S \subset M$ be an orientable\footnote{or course, the index can be defined as well for
any lagrangian submanifold, but since we are working here with Bohr - Sommerfeld lagrangian submanifolds which must be orientable  we give here the light version of the definition}  lagrangian submanifold
of a simply connected symplectic manifold and $\ga \subset S$ --- a loop. Since
the ambient symplectic manifold is simply connected, then there exists a disc
$D \subset M$, whose boundary coincides with $\ga$:
$$
\partial D = \ga.
$$
Topologically the discs with boundaries at $S$ represent elements of the group $\pi_2(M, S; \Z)$ thus for our given loop $\ga$ the space of such discs is discretized by the corresoinding integer data.

Then for any (almost) complex structure $I$, compatible with $\om$, one can
compare two trivializations of the anticanonical line bundle $K^{-1}$, restricted to
$\ga$, namely:

(A) the canonical trivialization, given by the orthogonal projection of the top
polivector field on $S$, which is dual to the volume form of the corresponding
riemannian metric,

(B) a trivialization, taken  over the disc $D$ and then restricted to the boundary
of the disc.

 The comparison gives a map
 $$
 \phi: \ga \to U(1),
 $$
 and the degree of this map gives some numerical index, which is called the Maslov index.
 It depends on the homotopy class of the loop and on the homotopy class of the disc
 which are encoded by the corresponding element of the group $\pi_2(M, S, \Z)$ thus
 it is not a pure invariant of the class of the loop. But instead it defines a map 
 $$
 \mu: \pi_2(M, S, ; \Z) \to \Z
 $$
and the point is that the topological type of the map doesn't depend on the almost complex structure,
which defines the trivialization (A), and is invariant uder lagrangian deformation of
$S \subset M$. Thus one gets a symplectic invariant (see \cite{KarMas}).

To compare the Maslov index with the universal Maslov class in the particular situation, studied
in the present paper, let's translate the definition of the Maslov index to the language of connections.
Trivialization (A) above is substituted by the canonical trivial connection $A_0 \in \sA_h(K^{_1}|_S)$;
instead of trivialization (B) one takes a flat hermitian connection $A_f \in \sA_h (K^{-1}|_D)$
which always exists over $D$ and which is unique up to gauge transformations. Thus, restricting
the picture to a choosen loop $\ga \subset S$, one gets two trivial connections $A_0|_{\ga},
A_f|_{\ga} \in \sA_h(K^{-1}|_{\ga})$ which are related by a gauge transformation, and the degree
of the transformation gives the value of the Maslov index for pair $(\ga, D)$. Equivalently
it is expressed by the formula
$$
\mu ([\ga, D]) = \frac{1}{2\pi i} \int_{\ga} (A_f|_{\ga} - A_0|_{\ga}) \in \Z.
$$

Let's illustrate
the definition by the following

\begin{prop} Let $S \subset M$ be an orientable lagrangian submanifold of a simply connected
symplectic manifold and a pair $(\ga, D)$ such that $\ga \subset S, D \subset M$ and
$\partial D = \ga$ represents element $[\ga, D] \in \pi_2(M, S, \Z)$. Then its Maslov index
$\mu([\ga, D])$ is trivial if and only if for any compatible almost complex structure
$I$ there exists a flat hermitian connection $A_f \in \sA_h(K^{-1}|_D)$ such that
$$
A_f|_{\ga} = A_0|_{\ga},
$$
where $A_0$ is the canonical trivial connection, induced by $I$.
\end{prop}

In other words, the Maslov index is trivial if and only if the canonical trivial connection $A_0$
can be extended from $\ga$ to a flat connection on $D$.

Now come back to the case of Bohr - Sommerfeld submanifolds of a pseudo - Einstein
symplectic manifold. We have the following geometric interpretation of the universal Maslov class,
which looks quite similiar to the previous proposition about the Maslov index.

\begin{prop} Let $S \subset M$ be a Bohr - Sommerfeld lagrangian submanifold of a simply
connected compact pseudo - Einstein symplectic manifold with the universal Maslov class
$m_S \in H^1(S, \Z)$, and let $\ga \in S$ be a smooth loop in $S$ representing class
$[\ga] \in H_1(S, D)$. Then $m_S([\ga]) = 0$ if and only if for any disc $D \subset M,
\partial D = \ga$ and any compatible almost complex structure $I$ there exists a connection
$A \in \sA_h(K^{-1}|_D)$ such that $F_A = -2 \pi i k \om|_D$ and $A|_{\ga}
= A_0|_{\ga}$, where $A_0 \in \sA_h(K^{-1}|_S)$ is the canonical trivial connection
on $S$.
\end{prop}

In other words, a loop $\ga \subset S$ has trivial Maslov class value if and not only if
the loop can be glued by a disc, but if and only if
  the pair $(\ga, A_0|_{\ga})$ (which is often called {\it supercycle})
can be extended to a pair $(D, A \in \sA_h(K^{-1}|_D))$ such that $F_A = -2 \pi i k \om|_D$.

{\bf Proof.} Suppose $m_S([\ga]) = 0$. This means that for an almost complex structure $I$ there
exists some reference connection $a \in \sA_h(K^{-1})$ with curvature $F_a = - 2 \pi i k \om$
such that the difference
$$
a|_{\ga} - A_0|_{\ga} = 2 \pi i \rho
$$
is a pure imaginary closed 1 - form with trivial integral over $\ga$:
$$
\int_{\ga} \rho = 0
$$
(since it is exactly the value of the Maslov class on $[\ga]$). This implies
that 1 - form $\rho$ is exact and there exists a real smooth function $f$ such that
$$
d f = \rho.
$$
For any disc $D$, lies inside $\ga$, there exists an extension of function $f$ on
the boundary $\partial D = \ga$ which is a smooth function $F$ on whole disc $D$.
Then the required connection $A$ on $D$ is given by the formula
$$
A = a|_D - 2 \pi i d F \in \sA_h(K^{-1}|_D).
$$

In the opposite direction, let us fix any $I$ and $D$ and for this pair consider
the connection $A \in \sA_h(K^{-1}|_D)$ which exists by the statement. At the same time
let us take a reference connection $a \in \sA_h(K^{-1})$ over whole $M$. Since
the disc is simply connected and the curvatures of $a|D$ and $A$ are the same it follows
these connections are gauge equivalent each to other. Hence their difference
$$
\rho = \frac{1}{2 \pi i} (a|_D - A)
$$
is an exact real 1 -form. By the conditions of the statement $A|_{ga} =
A_0 |_{\ga}$ and consequently the value of the Maslov class $m_S$ on the cycle $[\ga]$
can be computed by the formula
$$
m_S([\ga]) = \int_{\ga} \rho = \int_{\ga} d f = 0,
$$
and it completes the proof.

Now one sees the difference between the Maslov index and the universal Maslov class:
while the first responses for the continuation to a flat connection over a given disc while the second
responses for the continuation to a connection with "harmonic" curvature form. The arguments,
already used in the proof of the last proposition, show that they are related in rather simple manner.

\begin{prop} In the situation, described above, for a given pair $(\ga, D)$ one has
$$
\mu ([\ga, D]) - m_S ([\ga]) = k \int_D \om.
\eqno (4)
$$
\end{prop}

The advantage of the Bohr - Sommerfeld condition is that for these lagrangian submanifolds
the symplectic area, multiplied by $k$, is an integer valued function on $\pi_2(M, S; \Z)$.
And if $k$ equals to zero, which happens in the case of trivial anticanonical class,
the universal Maslov class coincides with the Maslov index, and the last one doesn't depend
on the choice of $D$. 

To prove, let's take a compatible almost complex structure $I$ and consider three hermitian connections:
the canonical trivial connection $A_0 \in \sA_h(K^{-1}|_S)$, a flat connection $A_f \in \sA_h(K^{-1}|_D)$ and a reference connection $a \in \Oh_{k.\om} \subset \sA_h(K^{-1})$. Then
$$
 \int_D k . \om = \frac{1}{2\pi i}\int_D d (A_f - a|_D) = \frac{1}{2 \pi i} \int_{\ga}(A_f|_{\ga} -
A_0|_{\ga} + A_0|_{\ga} - a|_{\ga}) =  \mu ([\ga, D]) -  m_S([\ga]),
$$
and we are done.

Note that formula (4) looks quite similiar to the main relationship, derived in [2].
In some sense it is a generalization since we study non - integrable case, in some sense it is
a reduction since we claim it for Bohr - Sommerfeld lagrangian submanifolds only.

On the hand formula (4) hints another intepretation of the universal Maslov class. This interpretation
deals with monotone lagrangian submanifolds. The notion of monotonicity for
lagrangian submanifolds was introduced in \cite{Oh} and it plays important role in
the study of the Floer cohomology (see \cite{Oh}). Here we just remark that 

\begin{prop} A Bohr - Sommerfeld lagrangian submanifold of a compact simply connected
pseudo Einstein symplectic manifold is monotone if and only if its universal Maslov class
is trivial.
\end{prop}

It's clear from relation (4) that if the universal Maslov class $m_S$ vanishes then
 $S$ must be monotone.

Suppose now that a Bohr - Sommerfeld lagrangian submanifold is monotone.
It implies that for each pair $(\ga, D)$ the universal Maslov class of $[\ga]$ is proportional
to the symplectic area of $D$:
$$
 m_S([\ga]) = (\la - k) \int_D \om,
 \eqno (5)
 $$
 where $\la$ is the monotonicity coefficient, and it must hold for any $D$ for a fixed $\ga$.
 The contradiction comes with the fact that the left hand side doesn't depends on $D$.
 This means that for a fixed $\ga$ the symplectic area of all discs with the same boundary
 $\ga$ must be the same. But in the compact case it is not hard to find two different $D$s
 with different symplectic areas, and it is possible due to the fact that any
 compact symplectic manifold has nontrivial second cohomology and that the symplectic
 cohomology class is positive. Therefore the equality (5) is possible if and only if $\la = k$
 and $m_S([\ga])$ is trivial. Since it hold for any $\ga$, the Maslov class is trivial. 
 
 Thus in the compact case the universal Maslov class is an obstruction to monotonicity
 of Bohr - Sommerfeld lagrangian submanifolds.

At the same time one has
\begin{cor} If $S$ is a simply connected orientable lagrangian submanifold of a compact simply connected pseudo Einstein symplectic manifold then $S$ is always monotone.
\end{cor}

Indeed, every simply connected orientable lagrangian submanifold is Bohr - Sommerfeld, and
since $\pi_1(S)$ is trivial the universal Maslov class is trivial too.

\section*{Conclusion}

Remark 2 of Section 2 hints that there is a new relationship
on the space of Bohr - Sommerfeld lagrangian submanifolds of
a simply connected pseudo - Einstein manifold. This relationship
is defined as follows. Consider a pair of Bohr - Sommerfeld lagrangian
submanifolds $S_1, S_2$ with trivial universal Maslov classes:
$$
m_{S_1} = m_{S_2} = 0.
$$
Let us fix an almost complex structure $I$, compatible with $\om$.
Then it defines the corresponding pair of trivial connections
$A_0^{i} \in \sA_h(K^{-1}|_{S_i}), i = 1 , 2$. Since the universal Maslov classes
are trivial, each $A_0^i$ can be extended to a global reference connection $a^i
\in \sA_h(K^{-1})$. We say that $S_1$ and $S_2$ are equivalent
with respect to $(K^{-1}, I)$ if there exists a connection $a \in \sA_h(K^{-1})$
which is an extension of both $A_0^1$ and $A_0^2$:
$$
a|_{S_i} = A_0^i.
$$
Several questions arise: first, is the relationship is an equivalence indeed;
second, does this relationship depend on the choice of an almost complex structure
or it is universal; third, is this relationship a reduction of the standard
homology theory. One expects that these preliminary questions have meaningful
answers and this will lead to new interesting constructions. This theme will be in the focus
of our further investigations.

\end{document}